\documentclass[12pt]{article}

\usepackage[utf8]{inputenc}
\usepackage{graphicx,amsmath,amssymb,amscd,psfrag,stmaryrd}
\usepackage[sort&compress]{natbib}
\usepackage{dashrule}
\usepackage{authblk}
\usepackage{ifthen}

\newboolean{areva}
\setboolean{areva}{false}

\newcommand\ifa[2]
{\ifthenelse{\boolean{areva}}{#1}{#2}
}

\usepackage{tikz}
\newcommand{\fig}{}
\newcommand{\Cross}{\mathbin{\tikz [x=1.8ex,y=1.
8ex,line width=.1ex] \draw (0,0) -- (1,1) (0,1) -- (1,0);}}%
\DeclareMathOperator*{\cart}{\Cross}
\author{N. Desassis, D. Renard, H. Beucher, S. Petiteau, X. Freulon }

\ifa{\author[2]{S. Petiteau}
\author[2]{X. Freulon}
\affil[2]{AREVA}
}

\newcommand{\indexh}{{x_j-x_i\simeq \textbf{h}}}

\newcommand{\tb}{\textbf}
\newcommand{\h}{\textbf{h}}
\title{A pairwise likelihood approach for the empirical estimation of the underlying variograms in the plurigaussian models}
\date{}

\newcommand{\un}{1\hspace{-1.6mm}1}

\begin{document}

\maketitle

\abstract{
\noindent The plurigaussian model is particularly suited to 
 describe categorical regionalized variables.
  Starting from a simple principle, the thresholding of one or several 
 Gaussian random fields (GRFs) to obtain categories, the plurigaussian model is well adapted for a wide range of situations.
 By acting on the form of the thresholding rule and/or the threshold values (which can vary along space) and the variograms of the underlying GRFs, 
 one can generate many spatial configurations for the categorical 
 variables. One difficulty is to choose variogram model for the underlying GRFs. Indeed, these latter are  hidden by the truncation and we only observe the simple and cross-variograms of the category indicators. In this paper, we propose a semiparametric method based on the pairwise likelihood to estimate the empirical variogram of the GRFs. It provides an exploratory tool in order to choose  a suitable  model for each GRF and later to estimate its parameters.
 We illustrate the efficiency of the method with a Monte-Carlo simulation study\ifa{ and we apply it on one a real case study from the mining industry}{}.
 
\noindent The method presented in this paper is implemented in the R package RGeostats.
}

\vspace{1cm}

\hspace{-1cm} \textbf{Keywords: } Plurigaussian models ; empirical variography ; semiparametric estimation ;  pairwise likelihood (PL) ; underlying Gaussian Random Fields (GRFs)

\section{Introduction}


 Regionalized categorical variables often appear in several scientific domains. For instance, in the earth sciences, some continuous soil properties (e.g the permeability, the grade of an element, ...)  can be better described by first categorizing the rock types into lithofacies (or facies) which present a certain homogeneity with respect to the studied variable. Then, the continuous variables are studied separately in each category. In the scope of  conditional simulations, the lithofacies are first simulated conditionally to the observed lithofacies, then the continuous variables are simulated inside each simulated category according to their associated spatial distribution.
 To model and simulate a categorical random field, the  plurigaussian model is particularly appealing. Based upon the simple truncation of one  \citep{Mat1987} or several Gaussian Random Fields \citep{Lel1994,Lel1996}, it allows to reproduce a wide range of patterns. Applications of the plurigaussian model can be found for mineral resources evaluation  \citep{Car2007,Riq2008,Tal2015}, in hydrology \citep{Mar2009}. In petroleum, some authors use the plurigaussian models in link with history matching \citep{Hu2000,Liu2004,Rom2010}. In this paper, we suppose that we directly observe the categorical variable.

When the underlying Gaussian Random Fields (GRFs) are supposed to be stationary, two ingredients are necessary to fully specify the plurigaussian model: the coding function (or truncation rule) which defines the sets associated to each category and which can vary along the space and the multivariate covariance function.

Concerning the coding function, some authors use a simple parametric form, for instance a cartesian product of intervals and they allow the threshold values to vary along the space.
It is often the case vertically through the vertical proportion curves \citep{Fel2004} but also laterally, for instance when using an auxiliary information as seismic data in the model.
 Other authors concentrate on the estimation of more complex coding functions constant in space \citep{Ast2014}. Finally, \cite{All2012} estimate complex coding function varying in space by using auxiliary information.

 As mentionned by \cite{Mar2009}, one of the main difficulty arising from the use of the plurigaussian model is the inference of the variogram models of the underlying GRFs. Indeed, the  available empirical variograms are the variograms of the indicator functions of the categories (one simple variogram  per category  and  the cross-variograms for all the bivariate combinations) while the  variograms required by the model are the variograms of the underlying GRFs whose realizations are hidden by the truncation.

Until now, most of the methods to estimate the variogram of the underlying GRFs rely on the indicator variograms. For instance \cite{Mar2009} determine the variogram model of the underlying GRFs by using simulations. More precisely, they choose a parametric model for the underlying GRFs and they compute the parameters value such as the indicator simple variograms of the simulations are the closest to the data indicator variograms. The optimization is performed with simulated annealing. \citet{Arm2011}  exploit the mathematical relationships between the underlying GRFs variogram models, the coding function  and the indicator simple and cross-variograms. Some industrial softwares \citep[as][]{Isa2014} also use these relations and the users must choose the parameters of the variogram models of the underlying GRFs  by visual inspection of the resulting indicator variograms through  a trial-and-error procedure \citep{Gal1994}. \cite{Eme2007} performs the numerical integration of the Gaussian density by using its expansion into the normalized Hermite polynomials. 
All these methods are rather tedious as they have a high computational cost or require a lot of trials. \cite{Dow2003} and \cite{Xu2006} propose to find the range parameters automatically by minimizing a squared differences with a grid-search but the choice of the covariance models of the underlying GRFs remains arbitrary and limited.

In this paper, we will supppose that the coding function is known and we will concentrate on the estimation of the variograms of the underlying GRFs. We propose an original methodology based on the pairwise likelihood (PL) maximization principle to directly compute the empirical variograms of the underlying and hidden GRFs.  More precisely, we perform a semiparametric estimation by considering the variogram at a given lag (distance in the omnidirectional case or vector otherwise) as a parameter of the model. Then, we maximize the PL by grouping the pairs of points approximately separated by this lag in the same way as for the classical empirical variogram estimation \citep{Mat1962}. We iterate this calculation on all distances (respectively vectors). Thereby, we obtain an empirical variogram which helps the user to choose a suitable valid model that can then be fitted by least squares or estimated with a likelihood based approach. Then, the simple and cross-variograms in the indicator scale can be derived and compared to the empirical variograms of the indicators to check the quality of the resulting models.

In the first part, we give the main notations of the paper and we recall the definition of the plurigaussian model. In section \ref{vario}, the relationships between variograms of GRFs and variograms of indicators are recalled  for comparison purposes. Then we present our  method in section \ref{new}. First, we describe the general principle which should make the estimation of a complex multivariate spatial model possible. Then we describe with more details the implementation in the case where the underlying GRFs are considered as independent. To assess the  efficiency of the method and to evaluate the uncertainty associated to the variogram estimation, a Monte-Carlo study is performed and its results are summarized in section \ref{mc}.\ifa{The method is tested on a real case study from the mining industry in section \ref{cs}.}{} We finish with some perspectives offered by the method to simplify the inference of the hierarchical spatial models.

\section{The data model}

\subsection{General formulation of the plurigaussian model}

Let $\mathcal{F}=\{f_1,\dots,f_K\}$ be a finite set with $K$ categories.
For a set of $n$ sites $\{x_i\}_{1\leq i \leq n}$ of a domain $\mathcal{D}\subset\mathbb{R}^d$, we observe $\textbf{f}=(f(x_1),\dots,f(x_n))$  a $\mathcal{F}$-valued vector.
We suppose that for a given location $x\in\mathcal{D}$, the value $f(x)$ is the realization of a $\mathcal{F}$-valued random variable $F(x)$.

To characterize the spatial distribution of $F(.)$, we use the plurigaussian model as described in \cite{Arm2011}. Let $$\mathbf{Y}(.)=\{\mathbf{Y}(x),x\in\mathcal{D}\}$$ be a  $q$-variate centered and standardized GRF on $\mathcal{D}$: for all $x\in\mathcal{D}$,  $\mathbf{Y}(x)=(Y_1(x),\dots,Y_q(x))$ is a random vector with $q$ components and for all $N\in\mathbb{N}^\star$ and for all $(x_1,\dots,x_N)\in\mathcal{D}^N$, the $N\times q$-vector $$(\mathbf{Y}(x_1),\dots,\mathbf{Y}(x_N))=(Y_1(x_1),\dots,Y_1(x_N),\dots,Y_q(x_1),\dots,Y_q(x_N))$$ is a standard Gaussian vector with $E[Y_r(x_i)]=0$ and $\text{Var}[Y_r(x_i)]=1$ for all $r\in [\![1,q]\!]$ and $i\in[\![1,N]\!]$.
In this paper, we  supose that $\mathbf{Y}(.)$ is a second-order stationary multivariate function, i.e there exists a matricial cross-covariance function $\tb{C}$ such as $\textrm{Cov}(Y_r(x),Y_s(x'))=\tb{C}_{rs}(x'-x)$  for $(r,s)\in[\![1,q]\!]^2$ \citep[see][for an introduction on multivariate spatial random fields]{Wac2003}.

 Let $\mathcal{C}$ be a coding function on $\mathcal{D}$ such as, for all $x\in\mathcal{D}$, $\mathcal{C}(x)=(\mathcal{C}_1(x),\dots,\mathcal{C}_K(x))$ where, for $k\in[\![1,K]\!]$,  the subsets $\mathcal{C}_k(x)$  form a (measurable) partition of $\mathbb{R}^q$.
The model is defined by the following equivalence 

\begin{equation}\label{model}
F(x)=f_k \textrm{ if and only if } \mathbf{Y}(x)\in\mathcal{C}_k(x).
\end{equation}

 
Note that the formulation given by (\ref{model}) provides a quite general class of models. Indeed, it also contains the models defined by 
$$F(x)=f_k \textrm{ if and only if } \varphi(\mathbf{Y}(x))\in\tilde{\mathcal{C}}_k(x)$$

for any surjective function $\varphi$ from $\mathbb{R}^q$ to any set $E$ where the sets $\tilde{\mathcal{C}}_k(x)$ for $[\![1,K]\!]$ form a partition of $E$. The subsets $\varphi^{-1}(\tilde{\mathcal{C}}_k(x))$ have to be some measurable sets of $\mathbb{R}^q$. This remark aims to highlight the fact that the marginal gaussianity of the random variables $Y_r(x)$ is arbitrary. Nevertheless, the Gaussian assumption is a convenient way to describe the spatial multivariate relationships of the underlying random fields. It also provides a multivariate random field easy to simulate \citep[see e.g][]{Lan2002}.

We will note $c(x)$ the set defined as: $$c(x)=\mathcal{C}_k(x)$$ where $k\in[\![1,K]\!]$ is the index of the category at location $x$. In other words, $f(x)=f_k$. 
In all the sequel, we will suppose that the classes $\mathcal{C}_k(x)$ are known.
 
\subsection{Indicators cross-variograms based methods}\label{vario}

As mentionned in the introduction, most of the methods to choose the covariance model of the underlying GRFs rely on the indicator simple and cross-variograms or covariances.  \cite{Arm2011} or \cite{Isa2014}  use the mathematical relationships between the simple and cross-covariances (or variograms) of the GRFs and the simple and cross-covariances (or variograms) of the indicators of each category. In the current paper, we only use these relationships to check the quality of the results given by our proposed method. We recall these relationships below. For  that purpose, we will note the random indicator function of the category $f_k\in\mathcal{F}$ as follows:

$$\un_{f_k}(x) =\left\{\begin{array}{cl} 1 & \textrm{ if }  F(x)=f_k\\
0 & \textrm{ otherwise} \end{array}\right.$$

and $1_{f_k}(x)$ the associated true value. 

\subsubsection{Variogram between two points}

For $(k,l)\in[\![1,K]\!]^2$, one can define the variograms between indicators of facies $k$ and $l$, between two locations $x$  and $x'$ of $\mathcal{D}$ by:

    $$\gamma_{kl}(x,x')=\frac{1}{2}E[(\un_{f_k}(x')-\un_{f_k}(x))(\un_{f_l}(x')-\un_{f_l}(x))]$$

When $k=l$, we have the simple variogram:

\begin{equation}\label{v1}
  \gamma_{kk}(x,x')=\frac{E[\un_{f_k}(x)]+E[\un_{f_k}(x')]}{2}-E[\un_{f_k}(x')\un_{f_k}(x)]
\end{equation}

When $k\neq l$, we have the cross-variogram:

\begin{equation}\label{v2}
  \gamma_{kl}(x,x')=-\frac{E[\un_{f_k}(x')\un_{f_l}(x)]+E[\un_{f_l}(x')\un_{f_k}(x)]}{2}
\end{equation}

We note $\Sigma_x$ and $\Sigma_{x,x'}$ the respective correlation matrices of the vectors $\textbf{Y}(x)$ and $(\textbf{Y}(x),\textbf{Y}(x'))$.
Furthermore,  $g^{}_\Sigma(\textbf{u})$ stands for the centered and standardized Gaussian density with correlation matrix $\Sigma$ computed for the vector $\textbf{u}$.

With these notations, we can establish the link between $\gamma_{kl}(x,x')$ and the  correlations between the underlying GRFs. Indeed, the expectation of the indicator of facies $k$ (which corresponds to its proportion at location $x$) is equal to:

\begin{equation}\label{v3}
  E[\un_{f_k}(x)]=\int_{\mathcal{C}_k(x)}g^{}_{\Sigma_x}(\textbf{u})d\textbf{u}
\end{equation}

and 

\begin{equation}\label{v4}
  E[\un_{f_k}(x)\un_{f_l}(x')]=\int_{\mathcal{C}_k(x)}\int_{\mathcal{C}_l(x')}g^{}_{\Sigma_{x,x'}}((\textbf{u},\textbf{v}))d\textbf{u}d\textbf{v}
\end{equation}

where each integration symbol represents an integration over a $q$-dimensional space. These integrations and all the others mentionned in the current paper are integral of the Gaussian probability density function. They can be computed numerically with the efficient algorithm proposed by \citet{Gen1992}.

Note that it is sometimes useful to work with the non-centered covariances 
$E[\un_{f_k}(x)\un_{f_l}(x')]$ which can be computed in the same way. Indeed, it has the advantage to capture asymetry in the model. 

When the GRFs are stationary and the coding function $\mathcal{C}$ is constant over $\mathcal{D}$, the variograms of all the involved random fields only depend on the lag between the points. Therefore, we can deduce the simple and cross-variograms of the indicator for a given lag from the variograms value of the underlying GRFs  by using formulas (\ref{v1}), (\ref{v2}), (\ref{v3}) and (\ref{v4}).
However, when the coding function varies over $\mathcal{D}$, the theoretical simple and cross-variograms of the indicators for a given lag do not exist anymore. Nevertheless, it is still possible to compute the  associated empirical variograms  and compare them with an averaged version of the variograms between two points computed in the indicators domain as described below.

\subsubsection{Variogram for a specific lag}

For a given vector $\h\in\mathbb{R}^d$, we will note $\indexh$  when the 
pair $(x_i,x_j)$ is used to compute the empirical variogram  at lag  $\h$ \citep[see e.g][for details]{Chi2012}. $N(\h)$ stands for the number of such pairs.

$$\hat\gamma_{kl}(\h)=\frac{1}{2N(\h)}\sum_\indexh(1_{f_k}(x_j)-1_{f_k}(x_i))(1_{f_k}(x_j)-1_{f_k}(x_i)).$$

This is the quantity that many authors suggest to fit with the image of the variogram model of the GRFs in the indicators scale which is defined for $(k,l)\in[\![1,K]\!]^2$ by:

\begin{equation}\label{varioexpindic}
  \gamma_{kl}(\h)=\frac{1}{N(\h)}\sum_\indexh\gamma_{kl}(x_i,x_j).
\end{equation}

This quantity depends on the spatial characteristics of $\mathbf{Y}(.)$  (defined through its multivariate cross-covariance function in the stationary case) and from the the set functions $\mathcal{C}_k$ through the values at the observation locations $\{x_1,\dots,x_n\}$ 

\section{Estimating the variogram by using pairwise likelihood}\label{new}

In this part, we describe a new methodology to perform the multivariate empirical variography of the underlying GRFs from the category observations.
This methodology is based on the pairwise likelihood (PL) maximization.
We first recall the principle of the more general composite likelihood based approach. Then we show how to apply it for the plurigaussian model. Finally, we describe more precisely the algorithm in two particular cases: the monogaussian case ($q=1$) and the plurigaussian case in which the $q$ GRFs are independent and the sets $c(x)$ are cartesian products of real subsets.

\subsection{Composite likelihood maximization}\label{cl}

The PL approach belongs to the family of the composite likelihood methods \citep[see e.g.][for a comprehensive review]{Var2011}. It is generally used to estimate a parameters vector $\theta$ of a statistical model, for instance when the usual maximization of the full likelihood is computationally cumbersome. In these cases, the full likelihood is replaced by a weighted product of marginal or conditional likelihoods.  \cite{Lin1988} defines the composite likelihood as follows: if $W$  a random vector with multivariate density $f(w;\theta)$ and $\{\mathcal{A}_\beta,\beta\in\mathcal{I}\}$ is a set of marginal or conditional events with associated likelihoods $\mathcal{L}_\beta(\theta;w)\propto f(w\in\mathcal{A}_\beta;\theta)$ for a finite set $\mathcal{I}$, the composite likelihood is the weighted product $$\mathcal{L}_C(\theta;w)=\prod_{\beta\in\mathcal{I}}\mathcal{L}_\beta(\theta;w)^{\lambda_\beta}$$ where 
$\lambda_\beta$ are nonnegative weights to be chosen.

One of the advantages of the composite likelihood based approaches is that they enable to estimate only some components of the parameters vector $\theta$. We use this idea below to derive a semiparametric estimator of the underlying variograms. Before that, we present an introductive example to show that, under suitable assumptions,  the usual empirical variogram can be considered as a maximum of a composite likelihood.

\subsection{The empirical variogram as a solution of an optimization problem}

Suppose that we observe $(y(x_1),\dots,y(x_n))$  derived from an intrinsec GRF with variogram $\gamma$. 

\cite{Cur1999} model $\gamma$ with a parametric form $\gamma_\theta$ and estimate $\theta$ by maximizing the marginal likelihood based on pairwise differences:

\begin{equation}\label{mlpd}
  \mathcal{L}_C(\theta;w)=\prod_{i=1}^{n-1}\prod_{j=i+1}^{n}f(v_{ij};\theta)^{\lambda_{ij}}
\end{equation}

where $f(v_{ij};\theta)$ stands for the density of the increment $V_{ij}=Y(x_j)-Y(x_i)$ and $\lambda_{ij}$ are some weights to choose.

 Here we adopt a semiparametric approach, \citep[as][with the spectral density]{Im2007}. We suppose that $\theta$ contains the variogram values $\gamma(\h_\alpha)$ for a set of $n_l$ lags $\h_\alpha$, $\alpha\in[\![1,n_l]\!]$. Then we group the pairs of points $(x_i,x_j)$  according to their distance as follows: when their exists $\alpha$ such as  $x_j-x_i\simeq \h_\alpha$, we set $\lambda_{ij}$ to 1 and we consider that $\gamma(x_j-x_i)=\gamma(\h_\alpha)$;  otherwise, if it is not possible to  associate a pair $(x_i,x_j)$ to a lag, its weight is set to 0. 
Then, it  is straightforward to show that the quantity which maximizes the associated marginal likelihood based on pairwise differences (\ref{mlpd}) is nothing but the traditional empirical variogram estimator of \cite{Mat1962}. 
We use the same idea to estimate the underlying variograms in the plurigaussian model.

\subsection{Adaptation to the plurigaussian case}\label{pl}

 We would like to estimate  $\Sigma(\textbf{h}_\alpha)$, the cross-covariance matrices of the vectors $$(Y_1(x),Y_1(x+\h_\alpha),\dots,Y_q(x),Y_q(x+\h_\alpha))$$ for a set of $n_l$ separation vectors $\textbf{h}_\alpha, \alpha\in[\![1,n_l]\!]$.

The particular composite likelihood which is adapted to this problem is the pairwise likelihood

$$\mathcal{L}_C(\theta;w)=\prod_{i=1}^{n-1}\prod_{j=i+1}^{n}f(w_i,w_j;\theta)^{\lambda_{ij}}$$ where $f(w_i,w_j;\theta)$ are the bivariate densities of $(W_i,W_j)$ for all $(i,j)\in[\![1,m]\!]^2$.

To estimate the covariance matrices $\Sigma(\h_\alpha)$, we consider that $\theta$ contains all the unknown elements of the matrices $\Sigma(\textbf{h}_\alpha)$ for $\alpha\in[\![1,n_l]\!]$. Then, we group the pairs of sites according to their separation vector in the same way as the empirical variogram computation and we write the log PL as follows:

\begin{equation}\label{pairwise}
  l_C(\theta;\textbf{f})=\sum_{\alpha=1}^{n_l}\sum_{\hspace{2mm}x_j-x_i\simeq\textbf{h}_\alpha}\textrm{log } p_{ij}(\Sigma(\textbf{h}_\alpha))
\end{equation}

where $$p_{ij}(\Sigma)=\int_{c(x_i)}\int_{c(x_j)}g^{}_\Sigma((\textbf{u},\textbf{v}))d\textbf{u}d\textbf{v}$$

is the probability that $F(x_i)=f(x_i)$ and $F(x_j)=f(x_j)$ when the cross-covariance matrix of the vector $(Y_1(x),Y_1(x+\h_\alpha),\dots,Y_q(x),Y_q(x+\h_\alpha)$  is $\Sigma$. Again, the weights $\lambda_{ij}$ attached to a pair $(i,j)$ have been set to 1 if there exists $\alpha\in[\![1,n_l]\!]$ such as $x_j-x_i\simeq \h_\alpha$ and to 0 otherwise.

Then, the PL estimator is obtained by maximizing $l_C$ with respect to all the matrices $\Sigma(\textbf{h}_\alpha)$. Note that to satisfy the stationarity of the resulting model, the condition $$\textrm{Cov}(Y_r(x_i),Y_s(x_i))=\textrm{Cov}(Y_r(x_j),Y_s(x_j))$$ is required for all locations $x_i$ and $x_j$ and all variable indices $(r,s)\in[\![1,q]\!]^2$. It implies that the $2q\times 2q$-matrices $\Sigma(\textbf{h}_\alpha)$ belong to the set noted $\mathcal{S}_{2q}$ defined by

 $$B\in\mathcal{S}_{2q}\Leftrightarrow b_{2r,2s}=b_{2r-1,2s-1}$$

 where $b_{r,s}$ stands for the $(r,s)^\textrm{th}$ element of the matrix $B$. 

 Furthermore, it is important to remark that all the matrices $\Sigma(\textbf{h}_\alpha)$  share some common terms to estimate, the ones corresponding to $\textbf{C}_{rs}(0)$. These two constraints on the global solution make the problem  numerically difficult to solve. 

In this paper, we  focus on the simplified cases where the $q$ GRFs are independent which is a quite ordinary assumption for the practitionners of the plurigaussian models. Then, $\textbf{C}_{rs}(\h)=0$ for all $(r,s)\in[\![1,q]\!]^2$ and all $\h\in\mathbb{R}^d$. Therefore, the matrices $\Sigma(\textbf{h}_\alpha)$ do not share any common term to estimate simultaneously. Furthermore, the positive definiteness of the matrices $\Sigma(\h_\alpha)$ is satisfied as soon as its non-zero off-diagonal elements belong to $]-1,1[$.
 
It results that the maximization of the log PL can be achieved by solving $n_l$ simpler maximization problems: for each $\alpha\in[\![1,n_l]\!]$, we can estimate the ``parameters''  $q$-vector $\theta_\alpha$ with $r^\textrm{th}$ element $\rho_r(\h_\alpha)=\textrm{Cor}(Y_r(x),Y_r(x+\h_\alpha))$ by

\begin{equation}\label{optim}
  \hat\theta_\alpha=\arg \!\!\!\max_{\hspace{-2mm}\rho\in]-1,1[^q}\sum_{\hspace{2mm}x_j-x_i\simeq\textbf{h}_\alpha} \textrm{log } \int_{c(x_i)}\int_{c(x_j)}\prod_{r=1}^qg^{}_{\rho_r}(u_r,v_r)du_rdv_r
\end{equation}
 
where $\rho=(\rho_1,\dots,\rho_q)$.

In the monogaussian case ($q=1$), the only quantity to estimate for a given lag $\h_\alpha$ is the spatial correlation of the underlying univariate GRF,  $\rho(\h_\alpha)=\textrm{Cor}(Y(x),Y(x+\h_\alpha))$ or equivalently $\gamma(\h_\alpha)=1-\rho(\h_\alpha)$.

Hence, the PL maximization problem (\ref{optim}) is reduced to a one dimensional optimization problem over a bounded interval; it can easily be solved, for instance with the golden section search algorithm \citep{Pre2007}.

For $q>1$, the generalization is straightforward if we assume that all the sets $\mathcal{C}_k(x)$ are cartesian products of subsets of $\mathbb{R}$:

$$\mathcal{C}_k(x)=\cart_{r=1}^q T_k^r(x)$$ with $T_k^r(x)\subset\mathbb{R}$, 

By denoting $t_r(x)=T_k^r(x)$ where $k$ is such that $f(x)=f_k$ is the actual category at site $x$, we have:

$$p_{ij}(\Sigma(\h_\alpha))=\prod_{r=1}^q\int_{t_r(x_i)}\int_{t_r(x_j)}g^{}_{\rho_r(\h_\alpha)}(u,v)dudv.$$

In other word, each $\rho_r(\h_\alpha)$ is estimated by:

$$\rho_r^\star(\h_\alpha)=\textrm{arg}\!\!\!\!\max_{\rho\in]-1,1[}\sum_{\hspace{2mm}x_j-x_i\simeq\h_\alpha} \textrm{log} \int_{t_r(x_i)}\int_{t_r(x_j)}g^{}_{\rho}(u,v)dudv.$$

which is equivalent to solve $q$ problems similar to the monogaussian case.

\section{Illustration on simulations}\label{mc}

In this section, we present two simulation studies to assess the efficiency of the proposed method. We work with two covariance models:

$$C_1(h)=e^{-h/20}$$
$$C_2(h)=e^{-(h/40)^2}$$

\subsection{The monogaussian case, $q=1$ }\label{mc1}

On a 1-dimensional regular grid with mesh size 1 and 2000 nodes, for $i=1,2$,  1000 realizations of standardized GRF $Y_i(.)$ with covariance  $C_i(h)$ have been drawn. 
 For each realization, $y(.)$, one category among the set $\mathcal{F}=\{\textrm{black},\textrm{red},\textrm{green}\}$ is assigned to each node $x$ of the grid according to the following rule:

   $$f(x)=\left\{\begin{array}{lll}
   \textrm{orange } &\textrm{ if } & y(x)\in\mathcal{C}_1(x)=(-\infty,s_1(x))\\
   \textrm{black }   &\textrm{ if } & y(x)\in\mathcal{C}_2(x)=(s_1(x),s_2(x))\\
   \textrm{green } &\textrm{ if } & y(x)\in\mathcal{C}_3(x)=(s_2(x),+\infty)
 \end{array}\right.$$

We consider two coding function cases:

 \begin{itemize}
 \item the constant case:  $s_1(x)=-s_2(x)$ are chosen such as the probability that $P(Y(x)\in\mathcal{C}_i(x))=\frac{1}{3}$ for all $i=1,2,3$ and all $x\in\mathbb{R}^d$, 
\item the varying  case where $s_1(x)$ and $s_2(x)$ have been simulated once for all simulations.
 \end{itemize}

Hence, 4 different categorical random fields have been considered by crossing all the possibilities ($C_1$ or $C_2$ vs. constant or varying coding function).
The empirical variogram of the underlying GRF is computed by pairwise likelihood from the categories as described in section \ref{pl}, for 150 distances ranging regularly from 1 to 150. For comparison purpose, in each case, the traditional empirical variogram has been computed directly on the realizations $y(.)$ for the same set of distances.

The results are summarized on figure \ref{result1}. In each case,  the average over all the simulations display a negligible bias. As expected, the variability of the estimator increases with the distance. The variogram seems to be better estimated when computed from categories by PL than with the original Gaussian values despite the loss of information due to the truncation. The reason is that we provide additional information by fixing the sill to 1 in the computation by PL. If we compare the results with respect to  the covariance model,  the simulations with model $C_1$ always display more important statistical fluctuations than with $C_2$. 
 The comparison between the constant coding function case and the varying one shows that the statistical fluctuations around the mean are greater for the latest. It is probably due to the fact that with the constant proportions $1/3$ of each category, a lot of transitions occur, bringing more information on the spatial correlation structure of the GRF than in the varying case. Indeed, for this latter case, the values of the probability of a given category are very strong in some areas, leading to few transitions and therefore less information on the hidden field.

\subsection{The independent bigaussian case, $q=2$ with $\mathcal{C}$  a constant cartesian product of intervals}\label{mc3}

In this example, we consider the same categories as previously. They are generated by assuming that $Y_1(.)$ and $Y_2(.)$ are independent with respective covariance function $C_1$ and $C_2$.
The categories are assigned to a point $x$ according to the following rule:

$$f(x)=\left\{\begin{array}{lll}
\textrm{black } &\textrm{ if } & y(x)\in\mathcal{C}_1(x)=(-\infty,s_1)\times\mathbb{R}\\
\textrm{orange }   &\textrm{ if } & y(x)\in\mathcal{C}_2(x)=(s_1,+\infty)\times (-\infty,t_1)\\
\textrm{green } &\textrm{ if } & y(x)\in\mathcal{C}_3(x)=(s_1,+\infty)\times (t_1,+\infty)
\end{array}\right.$$

where  $s_1=t_1=0$  such that 
$P((Y_1(x),Y_2(x)\in\mathcal{C}_1(x))=\frac{1}{2}$
and $P((Y_1(x),Y_2(x)\in\mathcal{C}_i(x))=\frac{1}{4}$   for $i=2,3$.

  A scheme of this coding function is displayed fig. \ref{flags} (b).

Then  1000 simulations are performed on 800 locations chosen uniformely on the square $[0,1]\times[0,1]$ one time for all the simulations. A realization is displayed fig. \ref{onerealbi}.

For each simulation, an empirical omnidirectional variogram is computed for 30 distances ranging from 0 and 150 with a tolerance factor on the distance. The results are summarized fig.\ref{result3} and again, they are rather good compared to the empirical variograms computed directly from the Gaussian data.
Note that for each simulation, the computation of the empirical variogram of the second Gaussian from $y(.)$ has been computed by using only the subset of locations for which the first Gaussian is greater than 0. 

\ifa{\section{Illustration on a mining example}\label{cs}

The case study concerns an uranium roll front deposit recongnized through 136 vertical drill holes 

 In this environment, the uranium deposit results from an oxidizing fluid flowing inside sandy deposit \citep{Lan2008}.  The rock types are summarized by using  3 lithofacies adapted to this specific mineralization: the area where the sand is oxidized (represented in orange), the area where the uranium concentration is upper than a given value (represented in black) and the area not yet oxidized, said reduced area (represented in green). The lithofacies are characterized every meters along the vertical which is  30 meters long.
We can see the arrangement of these 3 lithotypes  as a sequential approach because their relative location results from a sweeping of the medium, the uranium zone being between the oxidized and the reduced one. Then in the plurigaussian model, only one underlying GRF is needed (see fig. \ref{flags} (a)). Another interpretation is to consider that the  uranium area superimposes the background deposit composed of the oxidized and reduced areas, because it is obtained with other measurements, grade above a given value and not geological description. In that case, one can use two GRFs, the first one splitting the domain from geological analyses and the second one from grades (see fig. \ref{flags} (b)).

In this paper, the proportions of each lithofacies are considered as stationary horizontally but variable along the vertical axis. They have been computed by horizontally averaging the indicators of each lithofacies and then by convolving a standardized Gaussian kernel. The corresponding proportions show an increasing proportion of the mineralized type in center of the unit (see fig. \ref{flags} (c)).

The  empirical variograms of the underlying GRFs are computed by PL for the two hypothesis above mentionned and then fitted under the constraint that the variance of the model is equal to 1. Two directions are considered, one horizontal (omnidirectional in the horizontal plane) and one vertical. The results are displayed on fig. \ref{expstudy}.

Then, in both case, the image of the models in the indicator scale are computed by  using equation (\ref{varioexpindic}). The comparison with  the indicator variograms are displayed figures \ref{varioindvert} and \ref{varioindhoriz}.

}

\section{Discussion}

 In this paper, we propose to use the pairwise likelihood principle to estimate empirical variograms of the underlying GRFs in the plurigaussian model. The use of a composite likelihood based approach to provide an exploratory data analysis tool is apparently original, even if, as above-mentionned, the usual empirical variogram  can be seen as the solution of an optimization involving a composite likelihood.

Once the empirical variograms of the underlying GRFs have been computed, we can use them to choose a valid multivariate variogram model. They can be fitted by least squares, for instance by using the algorithm proposed in \citet{Des2013} or their parameters can be estimated by a likelihood based method. The likelihood will probably remain intractable since it involves at least an integral on $\mathbb{R}^{n}$ where $n$ is the number of samples. Again, a composite likelihood based approach should be used instead. Again, the PL seems well suited.

To conclude, note that the method presented in this paper is implemented in the R-package RGeostats \citep{Ren2014} in the function named {\it vario.pgs}. Some demonstration scripts are provided through a tutorial on the dedicated website.

Further researches will concentrate on the generalization of the approach presented in this paper to the case where no independence assumption is made between the underlying GRF. Indeed, more complex transitions between categories can be investigated, with more general multivariate spatial models \citep[see][]{Gal2006}. In that case, all the elements (except the diagonal) of the correlation matrices $\Sigma(\h_\alpha)$ must be estimated with the constraints mentionned in section \ref{pl}. This is computationally much more challenging. Another natural extension could be to adapt the local variogram kernel estimator proposed in \cite{Fou2014} to the plurigaussian context by using PL. Then, it would be possible to deduce a non-stationary model for the underlying GRFs, for instance with varying anisotropies.

Finally, the PL likelihood approach to compute empirical variograms seems to be a promising idea which could be applied to other similar context of hidden variable, or variable known after a transformation.

To cite some of them:
 
\begin{itemize}
\item compute the empirical variogram of the underlying GRF in the hierarchical geostatistical models \citep[see e.g][]{Dig1998}. Some authors have already proposed a way to compute empirical variogram of underlying random fields in hierarchical models:  \cite{Oli1993} treats the binomial case, \cite{Mon2006} the poisson case. However, these estimators are based on the method of moments and the distribution of the underlying random field is not specified. Thus, the underlying intensity can only be predicted by kriging but they can not be simulated. An approach based on the PL in a distribution based framework could be a good alternative;
\item perform the multivariate empirical variography of the underlying GRFs when one has to deal with a continuous variable vs. discrete variable \citep{Eme2009}, or even two discrete variables \citep{Ren2008, Cha2011};
\item compute the empirical variogram of a variable at punctual level when the observations are some regularizations with different supports.
\end{itemize}

\hspace{-1cm} \textbf{Acknowledgments: } We are very grateful to Christian Lantuejoul for the idea to directly compute the empirical variogram of the underlying GRF. We would also like to thank Geovariances and the School of Earth Science of the University of Queensland to have partly funded this research during the visit of the first author in Brisbane in 2012/2013.
 
\bibliographystyle{chicago}
\bibliography{bibli}

\newpage

\begin{figure}
  \begin{center}
 \includegraphics[width=4cm]{\fig 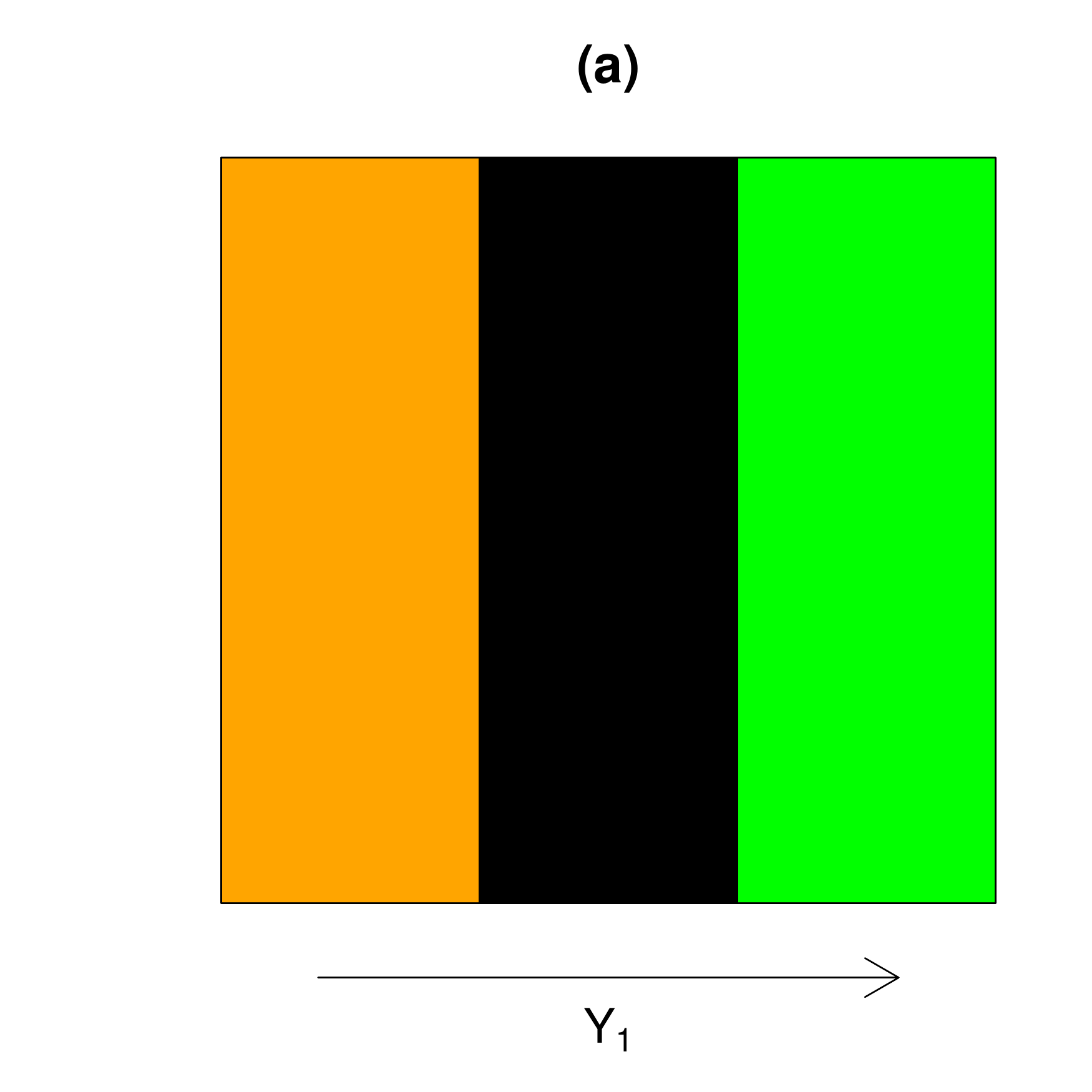}\includegraphics[width=4cm]{\fig 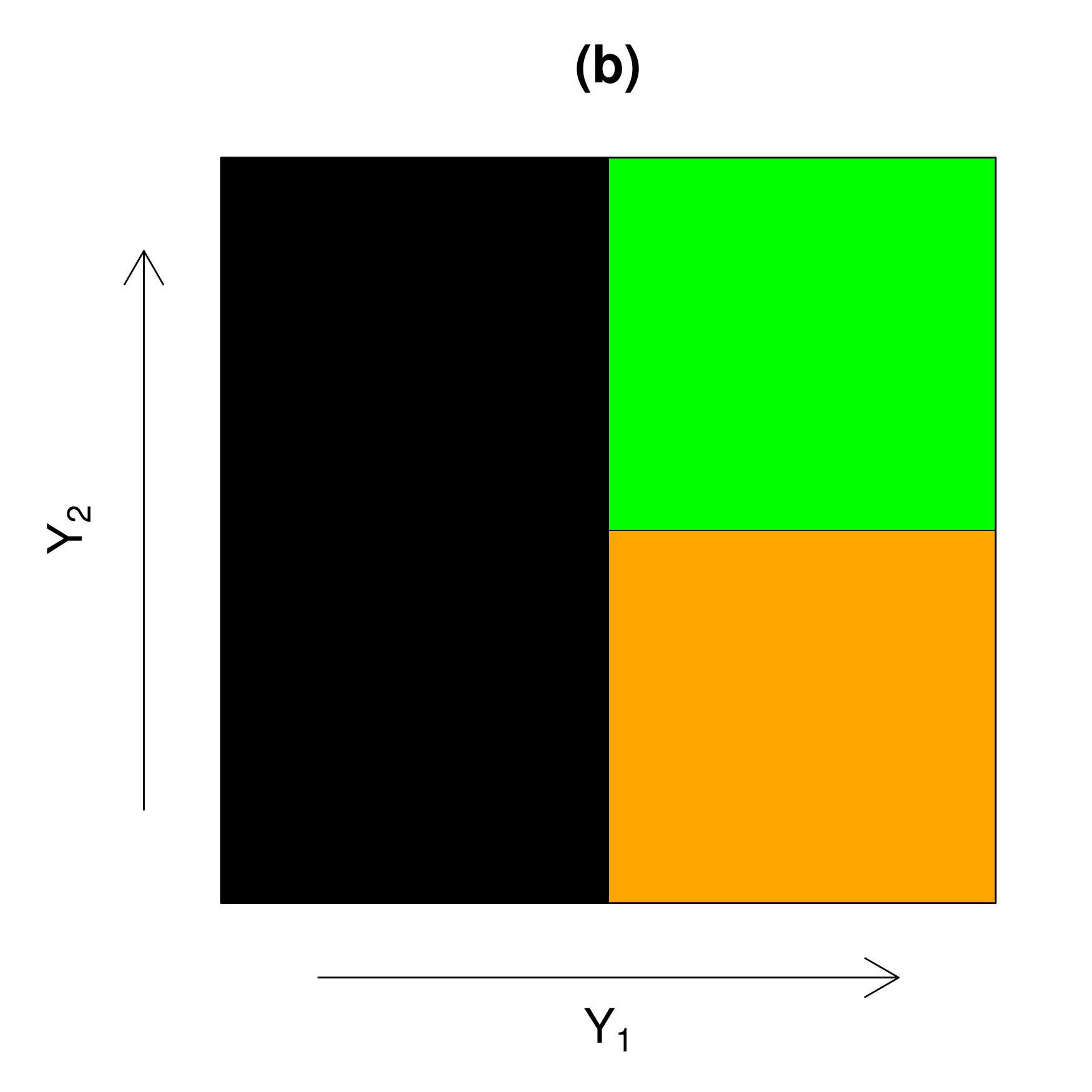}\ifa{\includegraphics[width=4cm]{\fig vpc.pdf}}{}
\caption[toc entry]{
Example of coding functions for (a) $q=1$ and (b) $q=2$. \ifa{(c) Vertical proportion curves for the case study.}

}\label{flags}
  \end{center}
\end{figure}

\newpage

\begin{figure}
  \begin{center}
   \includegraphics[width=3cm]{\fig 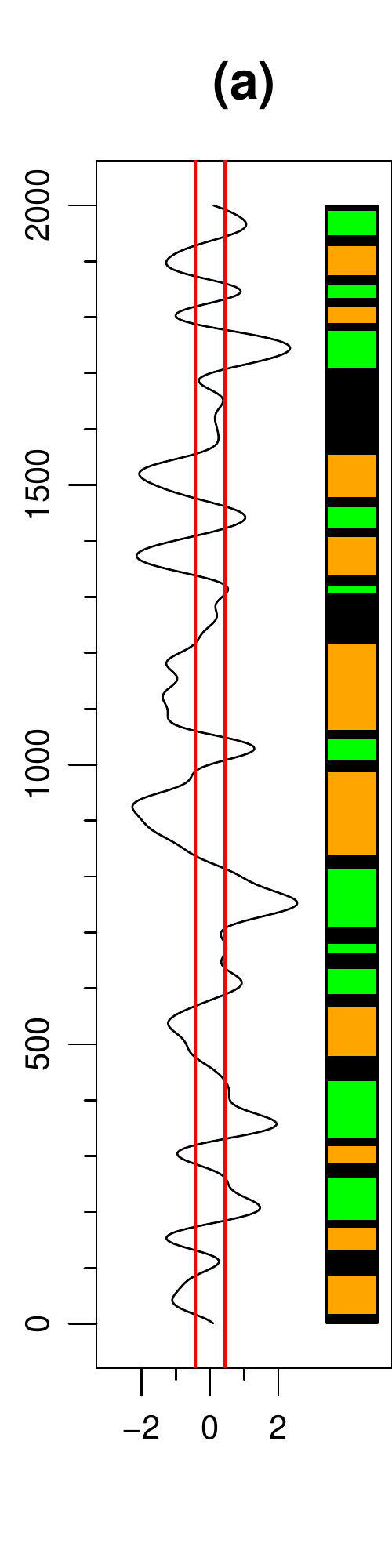}
 \includegraphics[width=3cm]{\fig 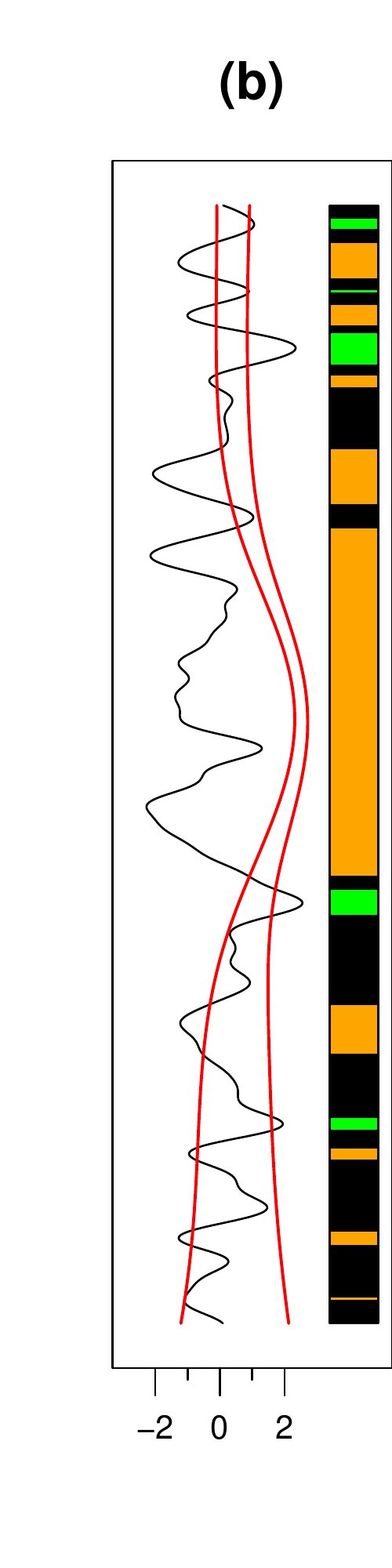}
\includegraphics[width=3cm]{\fig 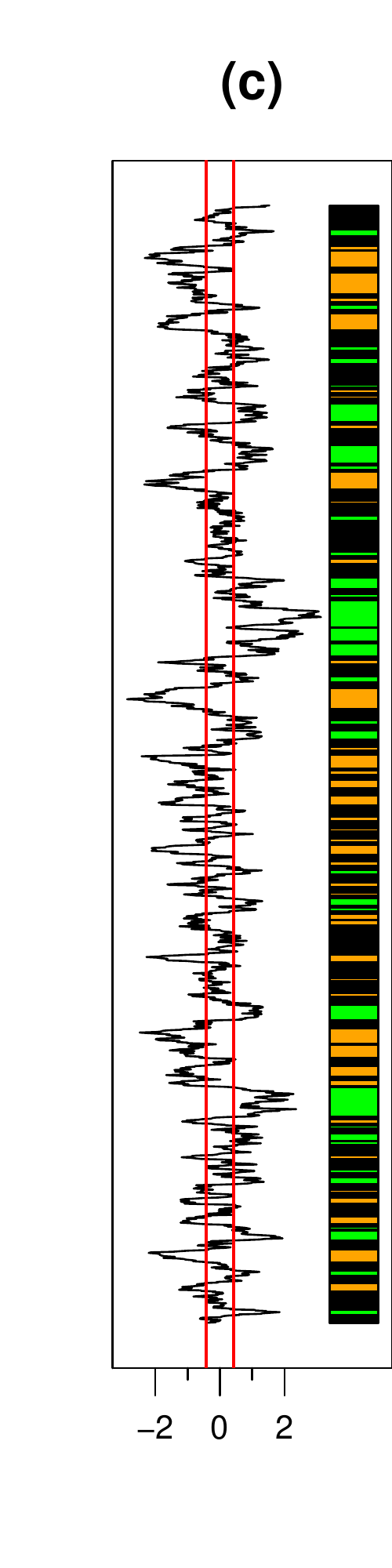}
\includegraphics[width=3cm]{\fig 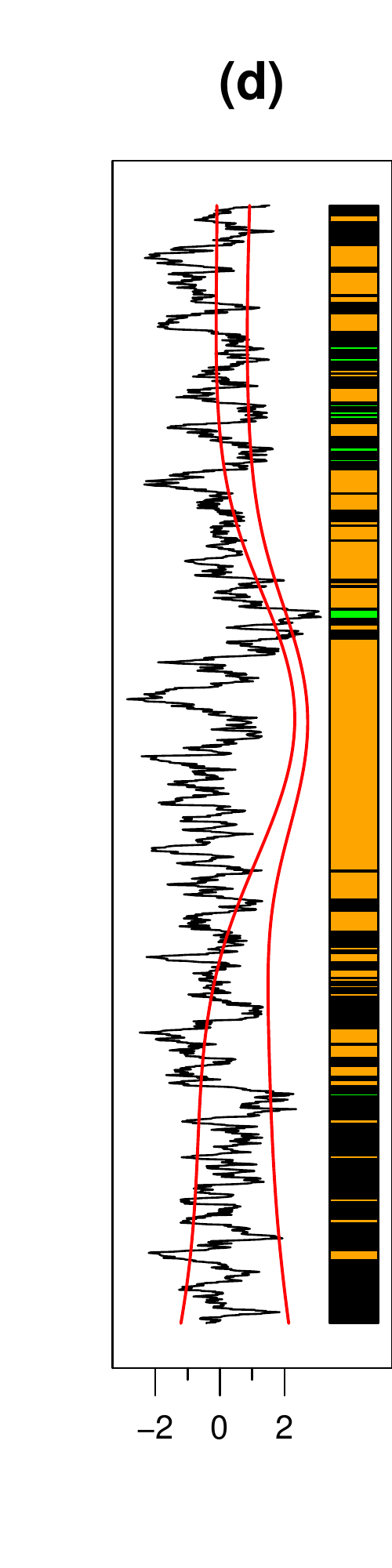}

\caption{For each configuration in the monogaussian case, one realization of $F$ and the associated $y$ (\hdashrule[0.5ex]{1cm}{0.5pt}{}), $s_1(x)$ and $s_2(x)$ ({\color{red}\hdashrule[0.5ex]{1cm}{0.5pt}{}}). The coding function is constant for (a) and (c), and variable for (b) and (d). The covariance model is $C_1$ for (a) and (b), and $C_2$ for (c) and (d).} 
\label{onereal}
  \end{center}
\end{figure}

\newpage

\begin{figure}
  \begin{center}
 \includegraphics[width=6.6cm]{\fig 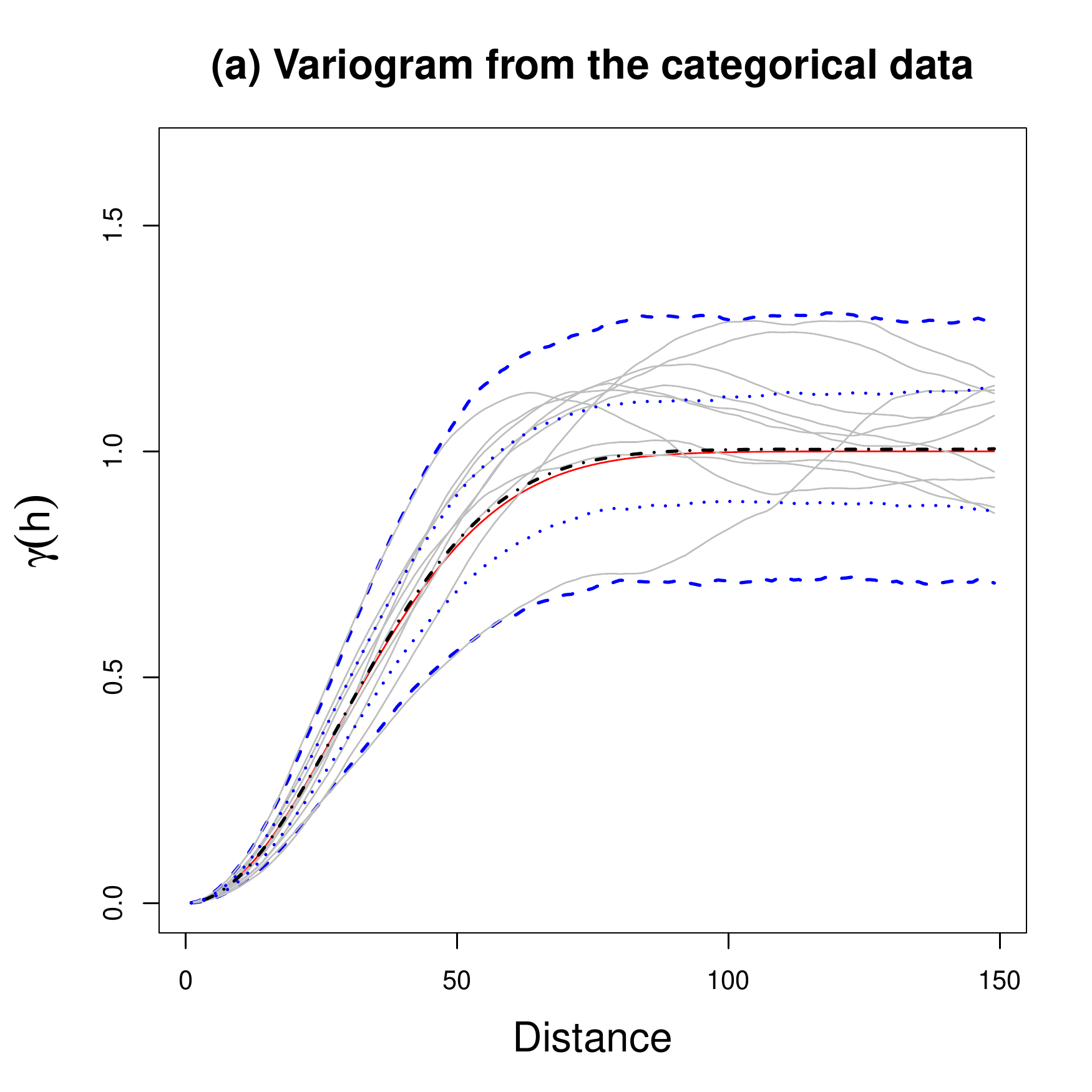}\includegraphics[width=6.6cm]{\fig 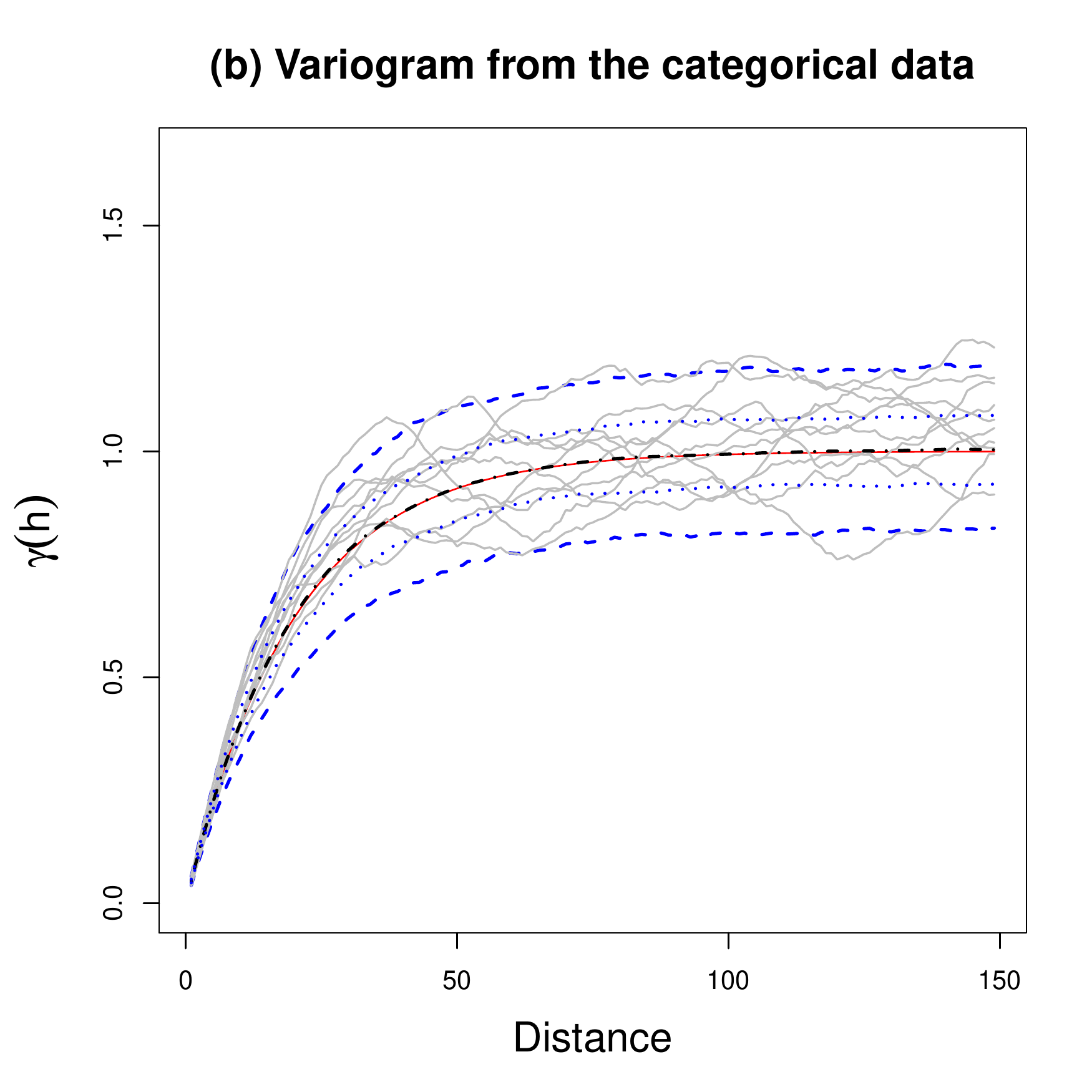}

\includegraphics[width=6.6cm]{\fig 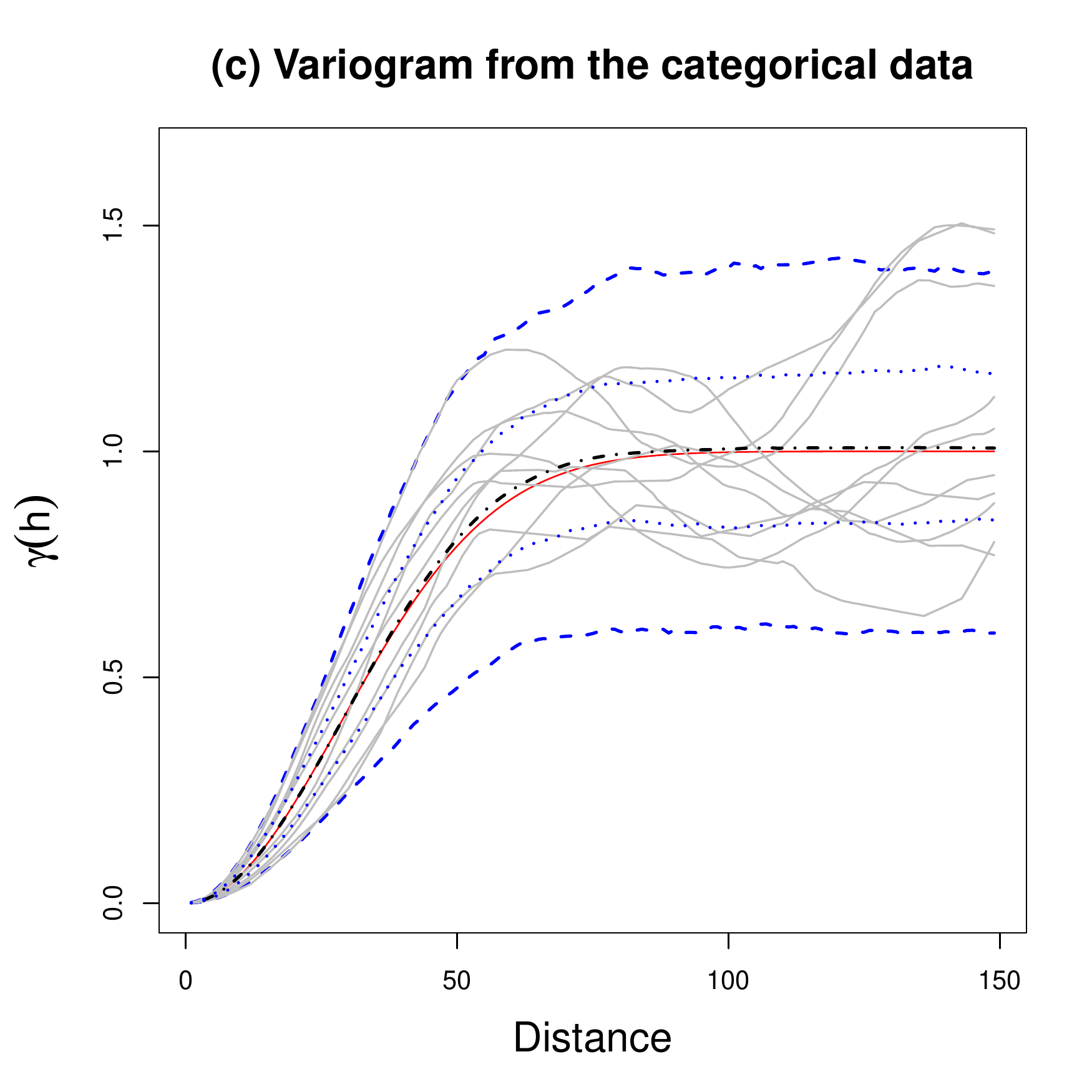}\includegraphics[width=6.6cm]{\fig 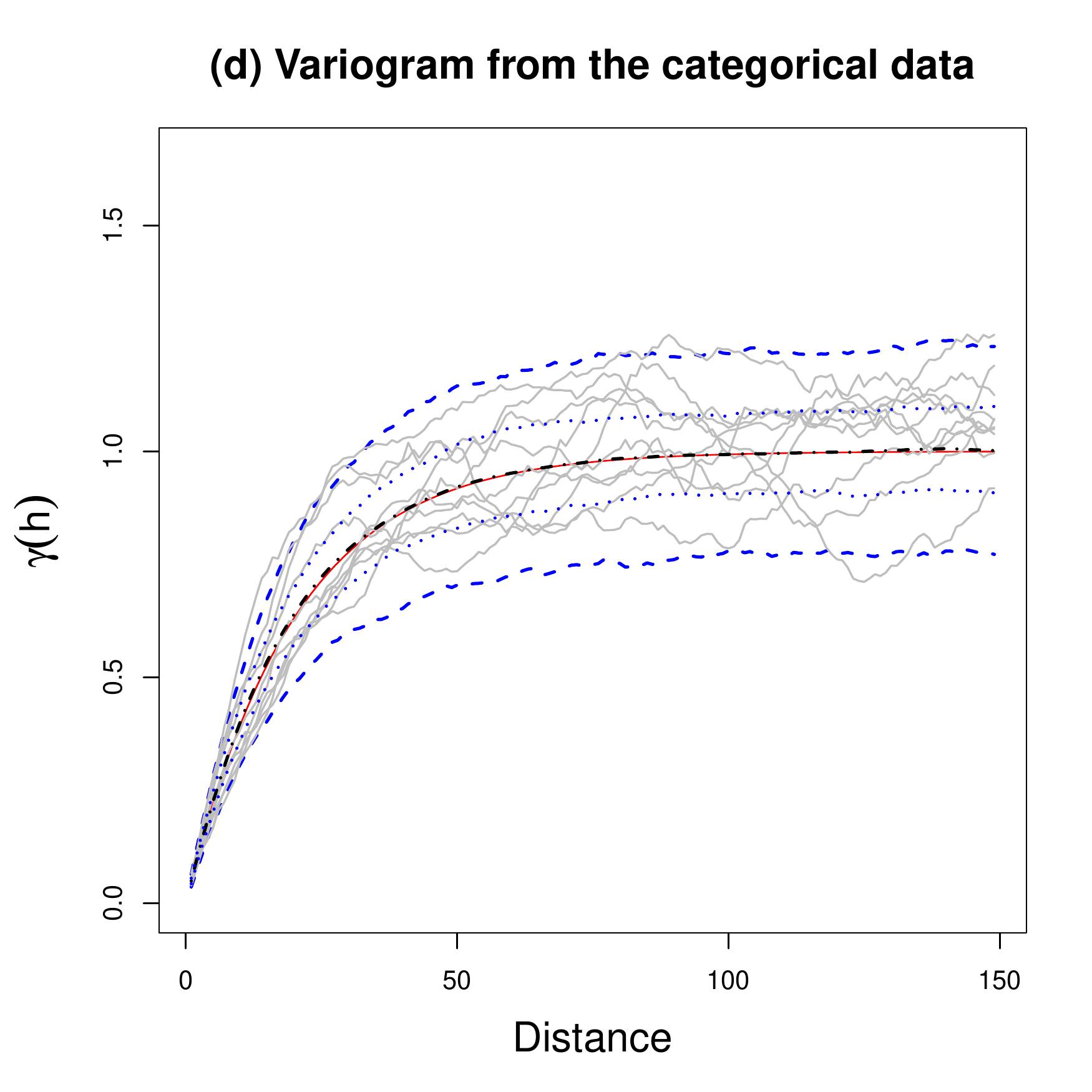}

\includegraphics[width=6.6cm]{\fig 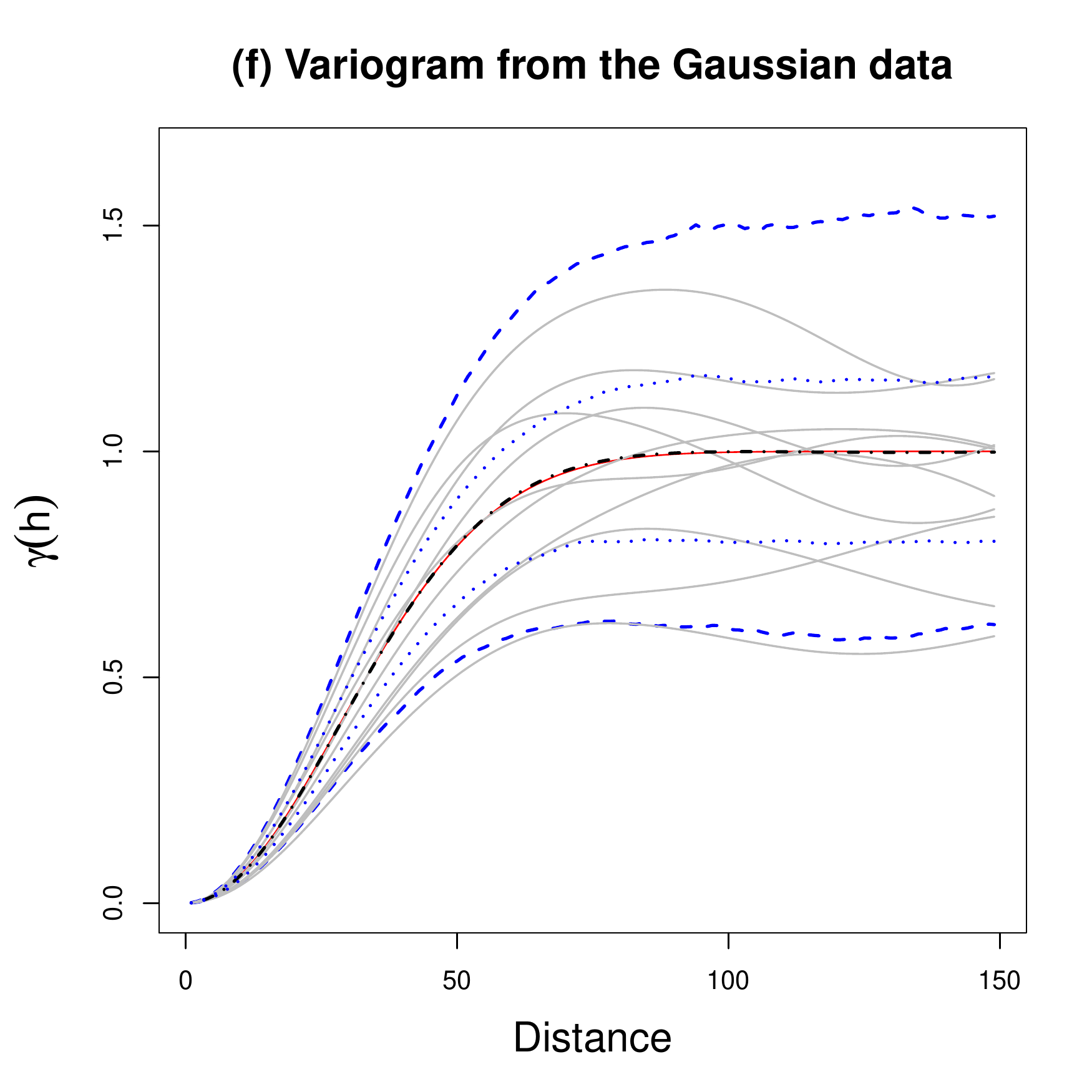}\includegraphics[width=6.6cm]{\fig 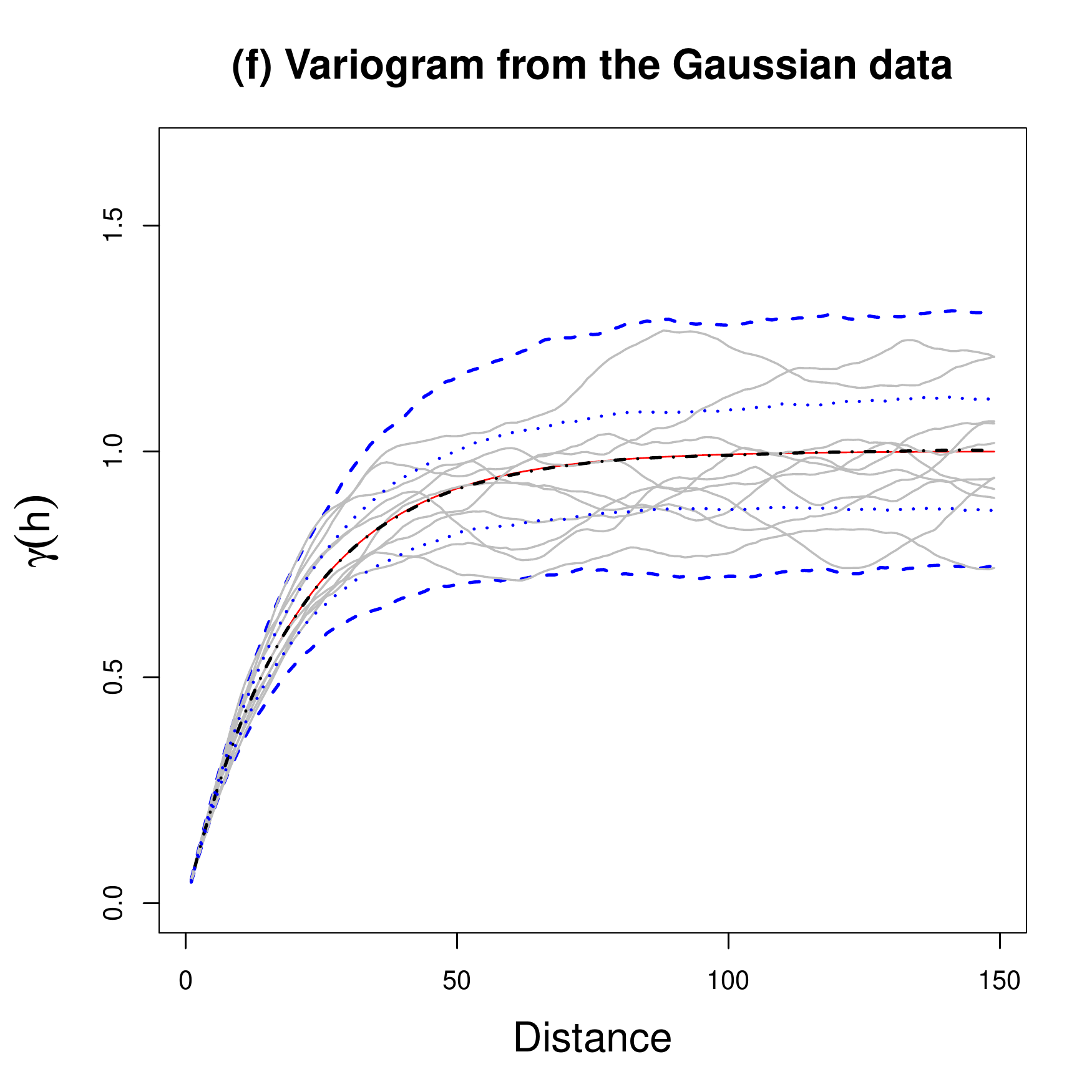}

  \caption{Simulation results in the monogaussian cases with covariance $C_1$ (left) and $C_2$ (right), in the constant coding function case (up) and in the varying case (middle). The down figures stand for the Gaussian case.
 Actual model ({\color{red}\hdashrule[0.5ex]{1cm}{0.5pt}{}}),
 empirical variogram of ten arbitrary simulations ({\color{gray}\hdashrule[0.5ex]{1cm}{0.1pt}{}}), 
average of the empirical variograms over all the simulations
(\hdashrule[0.5ex][x]{1cm}{1pt}{2.5mm 1mm 1pt 1mm}), 
 25th and 75th percentiles  ( {\color{blue}\hdashrule[0.5ex]{1.2cm}{1pt}{1pt 1mm}}),
 5th and 95th percentiles ( {\color{blue}\hdashrule[0.5ex][x]{1.2cm}{1pt}{1.5mm}}).}\label{result1}
  \end{center}
\end{figure}

\newpage

\begin{figure}
  \begin{center}
   \includegraphics[width=12cm]{\fig 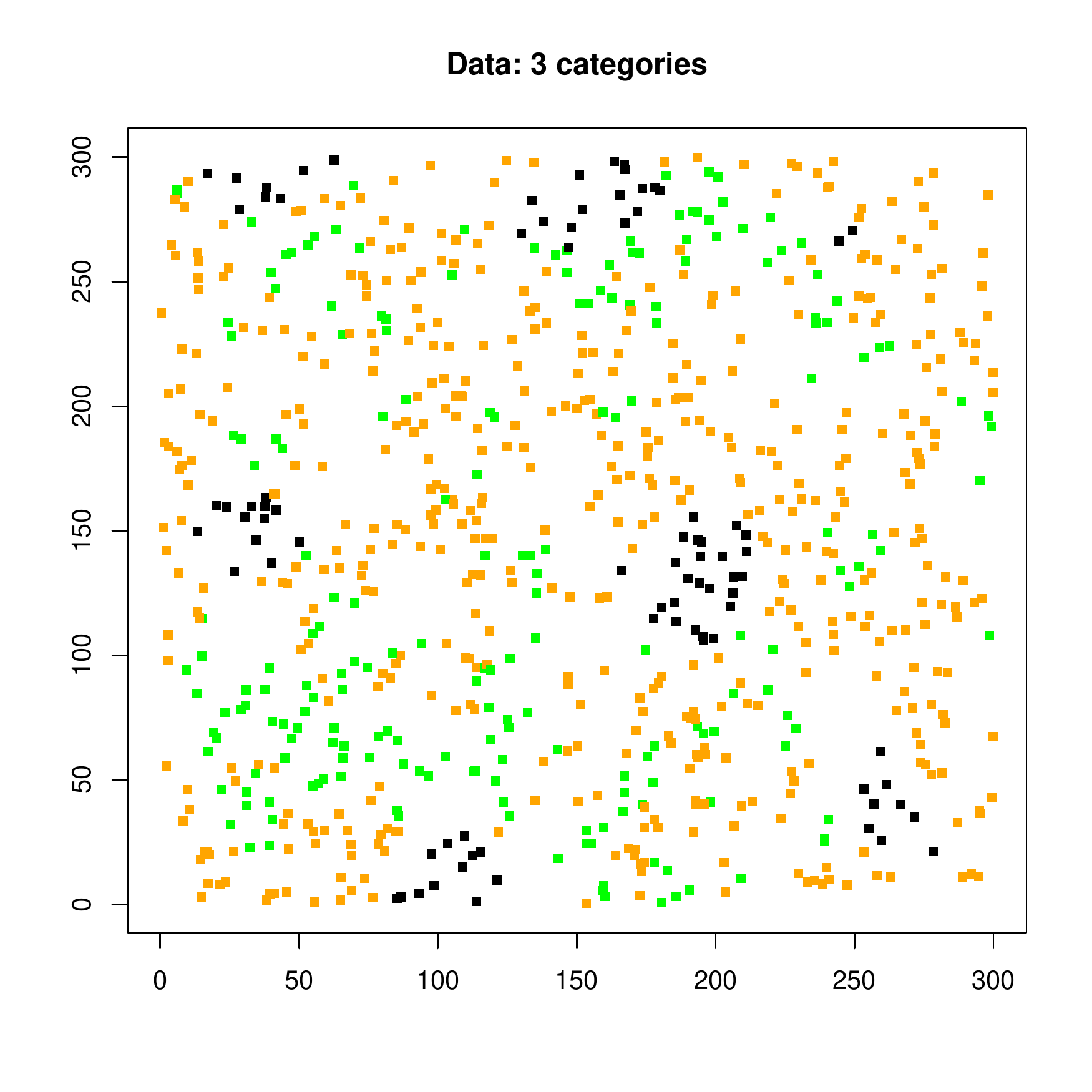}
  \label{onerealbi}
\caption{Location of the simulated points. The color indicates the simulated category for one specific simulation.}  
\end{center}
\end{figure}

\newpage

\begin{figure}
  \begin{center}
 \includegraphics[width=6.6cm]{\fig 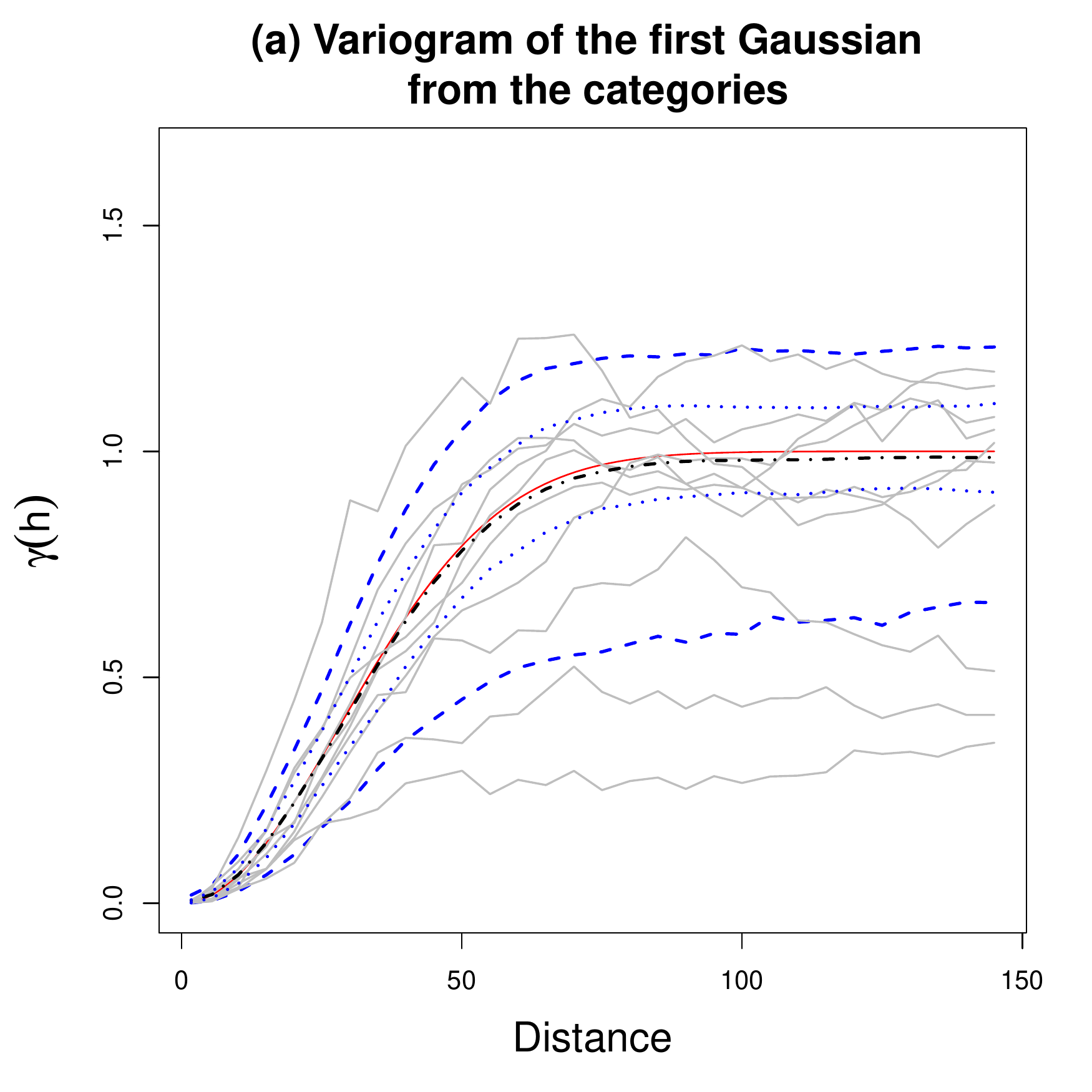}\includegraphics[width=6.6cm]{\fig 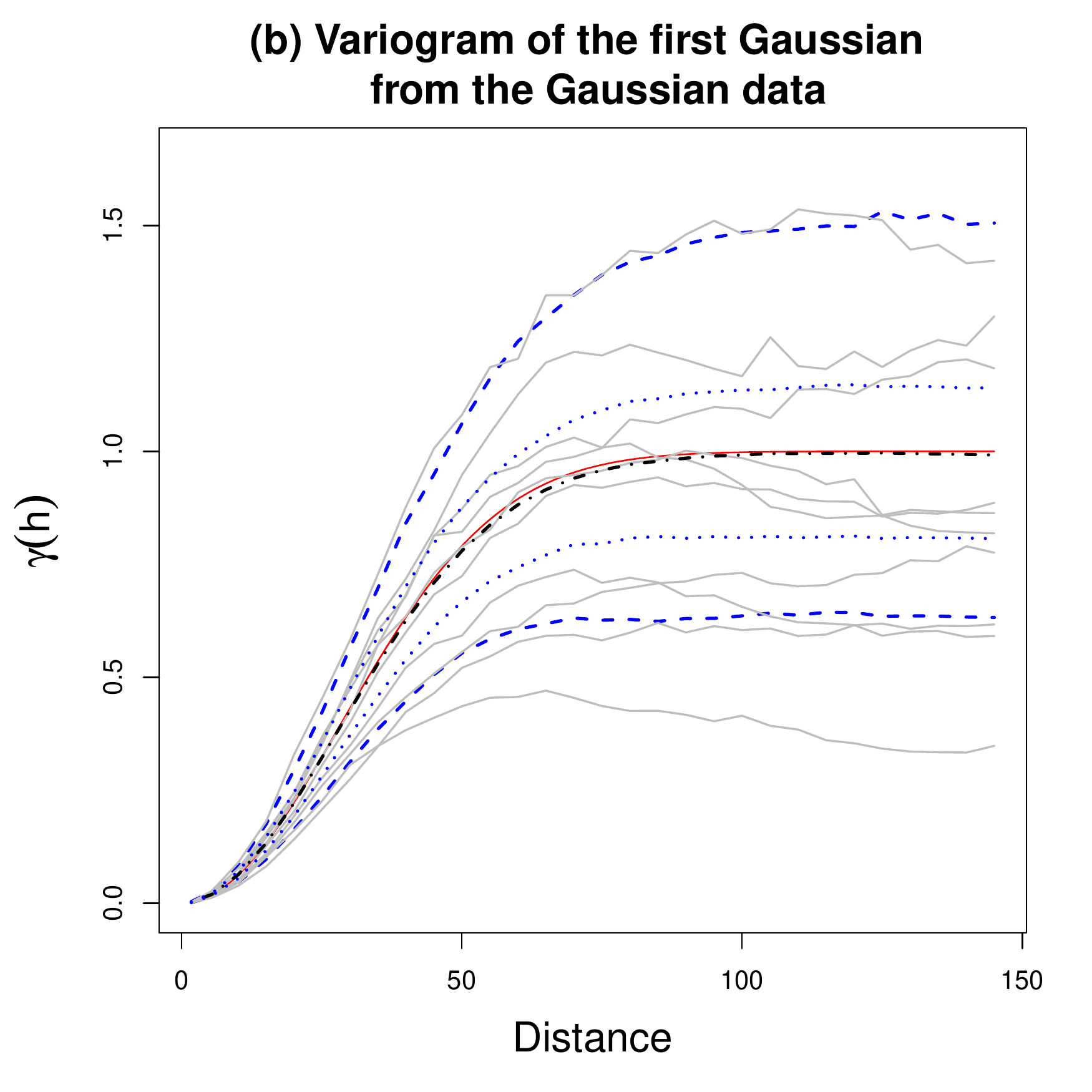}
\includegraphics[width=6.6cm]{\fig 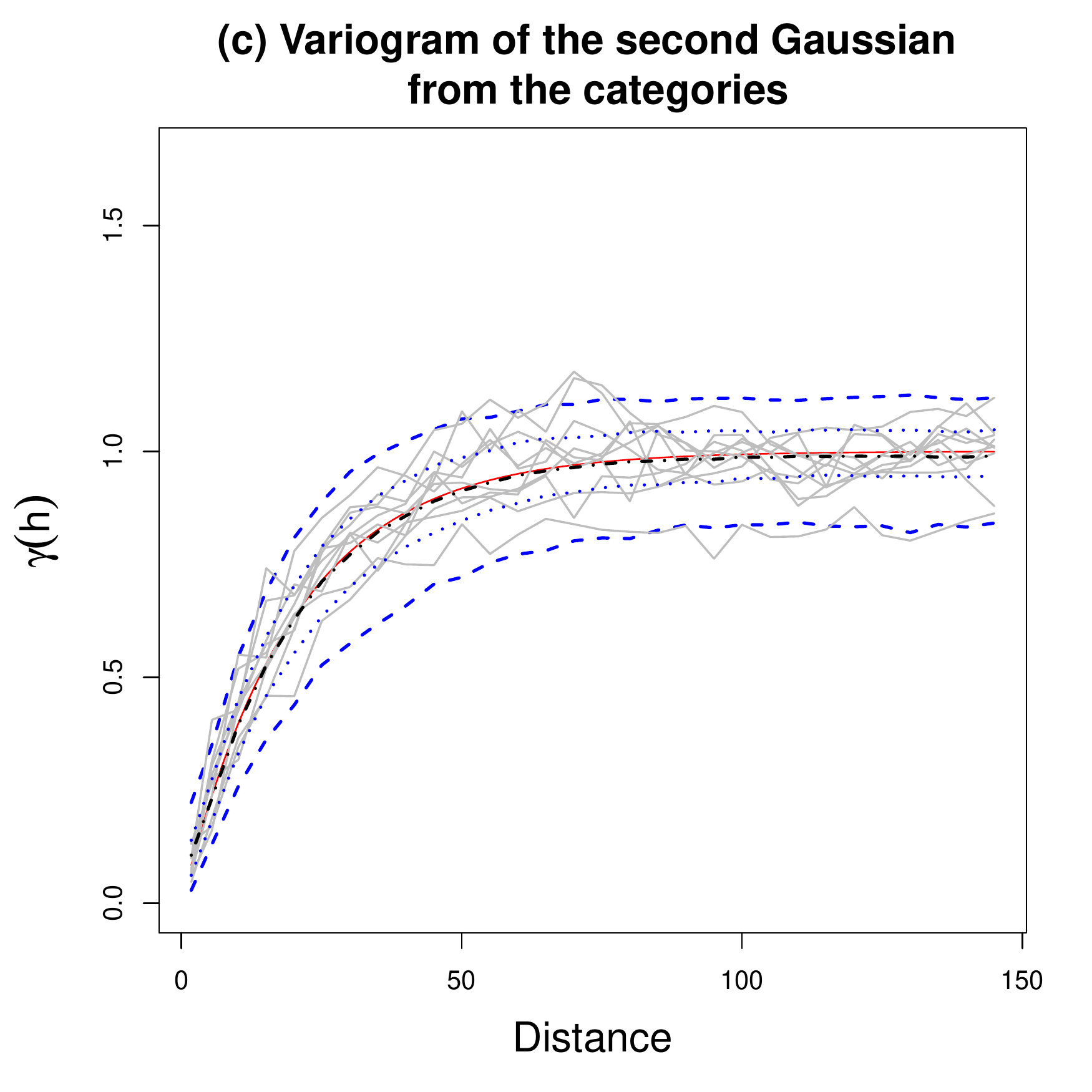}\includegraphics[width=6.6cm]{\fig 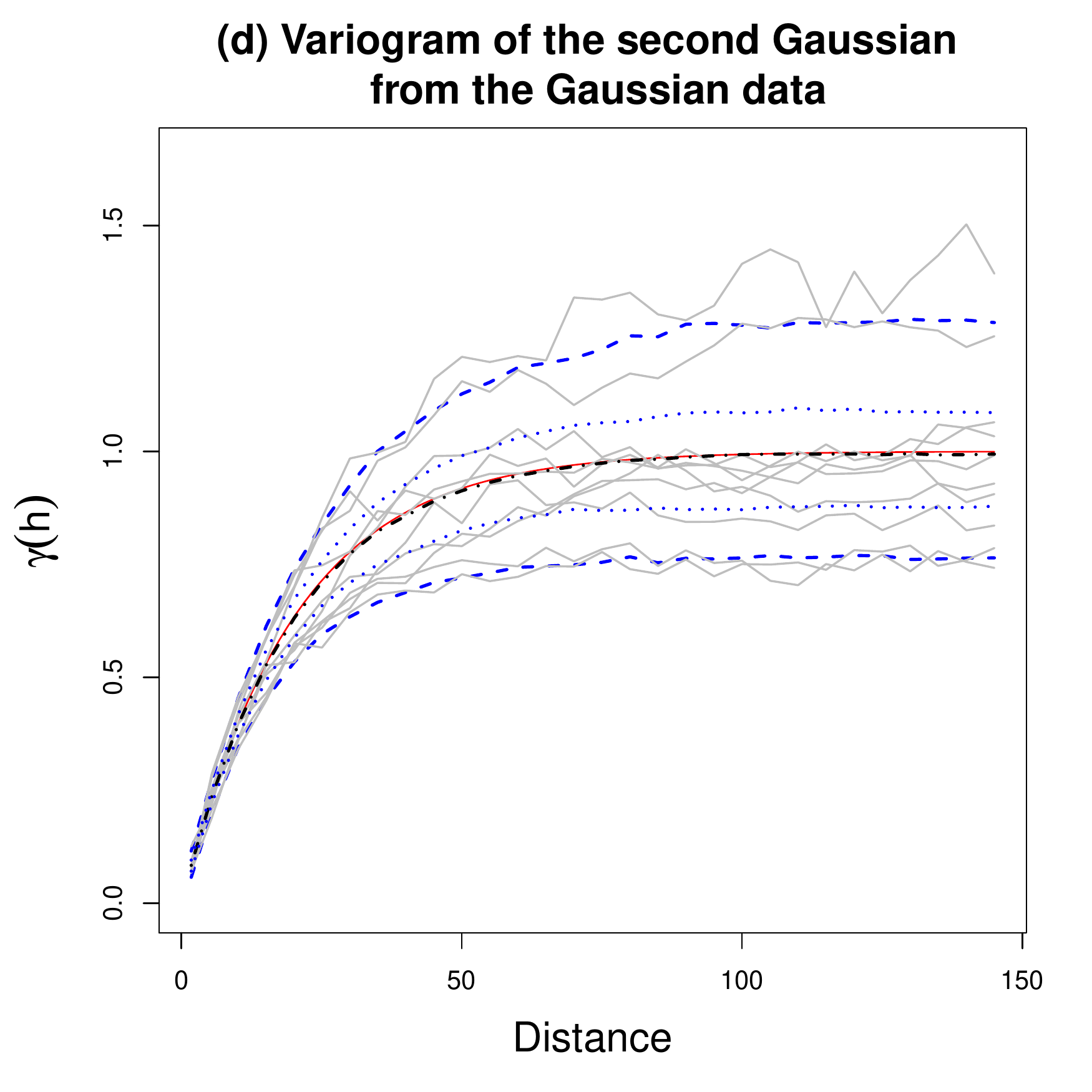}
  \caption{
 Actual model ({\color{red}\hdashrule[0.5ex]{1cm}{0.5pt}{}}),
 empirical ominidirectional variogram of ten arbitrary simulations ({\color{gray}\hdashrule[0.5ex]{1cm}{0.1pt}{}}), 
average of the empirical variograms over all the simulations
(\hdashrule[0.5ex][x]{1cm}{1pt}{2.5mm 1mm 1pt 1mm}), 
 25th and 75th percentiles  ( {\color{blue}\hdashrule[0.5ex]{1.2cm}{1pt}{1pt 1mm}}),
 5th and 95th percentiles ( {\color{blue}\hdashrule[0.5ex][x]{1.2cm}{1pt}{1.5mm}}).}\label{result3}
  \end{center}
\end{figure}




\ifa{}{\end{document}}

\newpage

\begin{figure}
  \begin{center}
\includegraphics[width=6.3cm]{\fig expgaussmonovert.pdf}\includegraphics[width=6.3cm]{\fig expgausspluri1vert.pdf}\includegraphics[width=6.3cm]{\fig expgausspluri2vert.pdf}

\includegraphics[width=6.3cm]{\fig expgaussmonohoriz.pdf}\includegraphics[width=6.3cm]{\fig expgausspluri1horiz.pdf}\includegraphics[width=6.3cm]{\fig expgausspluri2horiz.pdf}
  \caption{
Experimental variograms of the underlying GRFs obtained by PL and fitted model ({\color{red}\hdashrule[0.5ex]{1cm}{0.5pt}{}}) . Up : vertical; down: horizontal (omnidirectional in the horizontal plane).
From left to right: $q=1$, $q=2$ and first GRF, $q=2$ and second GRF.
}\label{expstudy}
  \end{center}
\end{figure}

\newpage

\begin{figure}
  \begin{center}
 \hspace{-9.7cm}\includegraphics[width=6.3cm]{\fig indicstudy11vert.pdf}

\hspace{-3.4cm}\includegraphics[width=6.3cm]{\fig indicstudy12vert.pdf}\includegraphics[width=6.3cm]{\fig indicstudy22vert.pdf}

\includegraphics[width=6.3cm]{\fig indicstudy13vert.pdf}\includegraphics[width=6.3cm]{\fig indicstudy23vert.pdf}\includegraphics[width=6.3cm]{\fig indicstudy33vert.pdf}

\caption{Vertical indicators simple and cross-variograms (from up to bottom and left to right: orange, black and green): computed by PL from the categorical data ({\color{black}\hdashrule[0.5ex]{1cm}{0.5pt}{}}), computed by equation (\ref{varioexpindic}) from the monogaussian model ({\color{green}\hdashrule[0.5ex]{1cm}{0.5pt}{}}), 
from the plurigaussian model
({\color{red}\hdashrule[0.5ex]{1cm}{0.5pt}{}}).}\label{varioindvert}
  \end{center}
\end{figure}

\newpage

\begin{figure}
  \begin{center}
 \hspace{-9.7cm}\includegraphics[width=6.3cm]{\fig indicstudy11horiz.pdf}

\hspace{-3.4cm}\includegraphics[width=6.3cm]{\fig indicstudy12horiz.pdf}\includegraphics[width=6.3cm]{\fig indicstudy22horiz.pdf}

\includegraphics[width=6.3cm]{\fig indicstudy13horiz.pdf}\includegraphics[width=6.3cm]{\fig indicstudy23horiz.pdf}\includegraphics[width=6.3cm]{\fig indicstudy33horiz.pdf}

  \caption{Horizontal indicators simple and cross-variograms (from up to bottom and left to right: orange, black and green) : computed by PL from the categorical data ({\color{black}\hdashrule[0.5ex]{1cm}{0.5pt}{}}), computed by equation (\ref{varioexpindic}) from the monogaussian model ({\color{green}\hdashrule[0.5ex]{1cm}{0.5pt}{}}), 
from the plurigaussian model
({\color{red}\hdashrule[0.5ex]{1cm}{0.5pt}{}}).}\label{varioindhoriz}

  \end{center}
\end{figure}

\end{document}